\documentclass[12pt]{amsart}
\usepackage[all]{xy}

\newtheorem{theorem}{Theorem}
\newtheorem{lemma}[theorem]{Lemma}

\newtheorem{cor}[theorem]{Corollary}

\newtheorem{conj}[theorem]{Conjecture}
\theoremstyle{definition}
\newtheorem{dfn}[theorem]{Definition}
\newtheorem{rk}[theorem]{Remark}
\newtheorem{ex}[theorem]{Example}
\def\<{\langle}
\def\>{\rangle}

\newcommand{\M}{\mathcal{M}}
\newcommand{\de}{\mathrm{d}}

\newcommand{\De}{\mathrm{D}}

\newcommand{\Z}{\mathbb{Z}}
\newcommand{\R}{\mathbb{R}}
\newcommand{\N}{\mathbb{N}}
\newcommand{\Q}{\mathbb{Q}}
\newcommand{\C}{\mathbb{C}}
\renewcommand{\to}{\rightarrow}
\newcommand{\To}{\longrightarrow}
\newcommand{\Mapsto}{\longmapsto}

\newcommand{\inclusion}{\hookrightarrow}
\newcommand{\p}{\partial}
\newcommand{\tf}{T_{\phi}}

\newcommand{\Of}{\Omega_{\phi}}
\newcommand{\ophi}{\omega_{\phi}}
\DeclareMathOperator{\id}{id}

\DeclareMathOperator{\im}{im}

\DeclareMathOperator{\ind}{index}

\DeclareMathOperator{\area}{area}

\DeclareMathOperator{\fix}{Fix}

\DeclareMathOperator{\symp}{Symp}
\DeclareMathOperator{\diff}{Diff}

\DeclareMathOperator{\sign}{sign}
\DeclareMathOperator{\flux}{Flux}
\DeclareMathOperator{\grow}{Growth}

\def\tr{\mathop{\rm tr}\nolimits}

\def\ind{\operatorname{ind}}

\newcommand{\HF}{HF_*}
\newcommand{\CF}{CF_*}

\def\Fix{\operatorname{Fix}}

\def\N{{\mathbb N}}

\begin{document}

\title[Floer homology and  Nielsen theory]
{Floer homology, Nielsen theory and symplectic zeta functions}

\author{Alexander Fel'shtyn}
\address{Fachbereich Mathematik, Emmy-Noether-Campus,
Universit\"at Siegen, Walter-Flex-Str. 3, D-57068 Siegen,
Germany}
\email{felshtyn@math.uni-siegen.de}
\urladdr{
http://www.math.uni-siegen.de/felshtyn}
\thanks{Part of this work was conducted during author stay in Ohio State
University, Columbus and  RIMS, Kyoto}
\begin{abstract}
We  describe
connection between symplectic  Floer homology for surfaces
 and Nielsen fixed point theory. A new  zeta functions and asymptotic invariant of symplectic
origin  are  defined. We show that   special values of symplectic
zeta functions  are  Reidemeister torsions.
\end{abstract}

\maketitle

\tableofcontents
\section{Introduction}
\providecommand{\bysame}{\leavevmode\hbox to3em{\hrulefill}\thinspace}
In two dimensions a diffeomorphism is symplectic if it preserves
area. As  a consequence, the symplectic geometry of surfaces lacks many of the interesting phenomena which are encountered in higher dimensions. For  example, two symplectic automorphisms of a closed surface are symplectically isotopic iff they are homotopic, by a theorem of Moser\cite{Mo}.
On other hand symplectic fixed point theory is very nontrivial in dimension 2,
as shown by the Poincare-Birkhoff theorem.
It is known that symplectic Floer homology on surface
is a simple model for the instanton Floer homology of the
  mapping torus of the surface diffeomorphism \cite{S}.

In this article we define new   zeta functions   related to symplectic Floer homology groups and investigate their analytical properties.  We show that  special values
of these zeta functions  are  Reidemeister torsions.
We also define asymptotic invariant of monotone symplectomorphism.
We  prove the following result: Let $M$ be a compact connected surface 
of Euler characteristic $\chi(M) <  0$.
If $\phi$ is  a non-trivial orientation
preserving periodic diffeomorphism of $M$ or  $\phi$ is  a  diffeomorphism of finite type   with  only isolated fixed points, then $\phi$ is monotone with
respect to some $\phi$-invariant area form on $M$ and
$
\HF(\phi) \cong
\Z_2^{N (\phi)},  \,\,\,  \dim\HF(\phi)= N (\phi)
$
where, $ N (\phi)$ denotes the Nielsen  number of $\phi$  
and $\HF(\phi)$ denotes symplectic Floer homology group.
Due to P. Seidel  \cite{S1}   $\dim\HF(\phi)$  is a new symplectic invariant
of a four-dimensional symplectic  manifold  with nonzero first Betti number.
This 4-manifold produced from symplectomorphism $ \phi$ by a surgery construction which is a variation of earlier constructions due to McMullen-Taubes,
Fintushel-Stern and J. Smith.We hope that our  symplectic zeta functions 
and asymptotic invariant  also
give rise to a new  invariants of contact 3- manifolds
and symplectic 4-manifolds. 

The author came to the  idea that Nilsen numbers are  connected with  Floer homology of surface diffeomorphisms  at the Autumn 2000, after conversations with Joel Robbin and Dan Burghelea. Last years I try  to find  a  cohomological theory which lies  behind of  the Nielsen zeta function \cite{f1}. The results of this paper were reported  in the  author talk on the International Conference Topological Methods in Nonlinear Analysis, June 2001 in Bendlewo, Poland.  I am very gratefull to Dietmar Salomon who send me handwritten notes of Wu-Chung Hsiang paper  ``A speculation on Floer theory and Nielsen theory'' and to  Stefan Haller, Wu-Chung Hsiang, Jarek Kedra,
Wilhelm Klingenberg,  Francois Laudenbach, Kaoru Ono,
Yuli Rudjak, Andrei Tyurin and  Vladimir Turaev  for very useful discussions of the results of this  paper. The delay in publication  was  connected with the  absence of the important  notion of monotonicity, introduced by P. Seidel\cite{S} in the  Spring 2001.

The author would like to thank the Max-Planck-Institute f\"ur Mathematik, Bonn
for kind hospitality and support.

\section{The monotonicity condition}
In this section we discuss the notion of monotonicity as defined in \cite{S,G}.
Monotonicity plays important role for  Floer homology in two
dimensions. For a more detailed account we refer to the original articles
of  P. Seidel and R. Gautschi.  

Throughout this article, $M$ denotes a closed connected and oriented
2-manifold of genus $\geq2$.   Pick an everywhere positive two-form $\omega$
on $M$. 

Let $\phi\in\symp(M,\omega)$, the group of  symplectic automorphisms 
 of the two-dimensional symplectic
manifold $M,\omega$.
The mapping torus of $\phi$,  $T_\phi = \R\times M/(t+1,x)\sim(t,\phi(x)),$
is a 3-manifold fibered over $S^1=\R/\Z$. 
There are two natural second
cohomology classes on $T_\phi$, denoted by $[\ophi]$ and $c_\phi$. The first one is
represented by the closed two-form $\ophi$ which is induced from the pullback
of $\omega$ to $\R\times M$. The second is the Euler class of the
vector bundle 
$$
V_\phi =
\R\times T M/(t+1,\xi_x)\sim(t,\de\phi_x\xi_x),
$$
which is of rank 2 and inherits an orientation from $TM$.

$\phi\in\symp(M,\omega)$ is called {\bf monotone}, if
$$
[\omega_\phi] = (\area_\omega(M)/\chi(M))\cdot c_\phi
$$
in $H^2(T_\phi;\R)$; throughout this article
$\symp^m(M,\omega)$ denotes the set of
monotone symplectomorphisms. 

Now $H^2(T_\phi;\R)$ fits into the following short exact sequence \cite{S,G}
\begin{equation}\label{eq:cohomology}
0 \To  \frac{H^1(M;\R)}{\im(\id-\phi^*)}
\stackrel{d}{\To} H^2(T_\phi;\R)
\stackrel{r^*}{\To} H^2(M;\R),
\To 0.
\end{equation}
where the map $r^*$  is restriction to the fiber. 
The map $d$ is defined as follows.
Let $\rho:I\to\R$ be a smooth function which vanishes near $0$ and $1$ and
satisfies $\int_0^1\!\rho\,\de t=1$.
If $\theta$ is a closed 1-form on $M$, then
$\rho\cdot\theta\wedge\de t$ defines a closed 2-form on $T_\phi$; indeed
$
d[\theta] = [\rho\cdot\theta\wedge\de t].
$
The map $r:M\inclusion T_\phi$ assigns to each $x\in M$ the
equivalence class of $(1/2,x)$.
Note, that $r^*\ophi=\omega$ and $r^*c_\phi$ is the Euler class of
$TM$. 
Hence, by \eqref{eq:cohomology}, there exists a unique class
$m(\phi)\in H^1(M;\R)/\im(\id-\phi^*)$ satisfying
$$
d\,m(\phi) = [\ophi]-(\area_\omega(M)/\chi(M))\cdot c_\phi,
$$
where $\chi$ denotes the Euler characteristic. 
Therefore, $\phi$ is monotone if and only if $m(\phi)=0$.
\smallskip\\
We recall the fundamental properties of $\symp^m(M,\omega)$ from \cite{S,G}. 
Let  $\diff^+(M)$ denotes  the group of orientation preserving
diffeomorphisms of $M$. 
\smallskip\\
(Identity) $id_M \in \symp^m(M,\omega)$.
\smallskip\\
(Naturality)\label{page:natur}
If $\phi\in\symp^m(M,\omega),\psi\in\diff^+(M)$, then 
$\psi^{-1}\phi\psi\in\symp^m(M,\psi^*\omega)$.
\smallskip\\
(Isotopy) 
Let $(\psi_t)_{t\in I}$ be an isotopy in $\symp(M,\omega)$, i.e. a smooth
path with $\psi_0=\id$.
Then 
$
m(\phi\circ\psi_1)=m(\phi)+[\flux(\psi_t)_{t\in I}]
$  
in $H^1(M;\R)/\im(\id-\phi^*)$; see \cite[Lemma 6]{S}.
For the definition of the flux homomorphism see \cite{MS}.
\smallskip\\
(Inclusion) 
The inclusion $\symp^m(M,\omega)\inclusion\diff^+(M)$ is a homotopy
equivalence.  
This follows from the isotopy property, surjectivity of the flux homomorphism
 and
Moser's isotopy theorem \cite{Mo} which says that each element of the mapping
class group admits representatives which preserve $\omega$.
Furthermore, the Earl-Eells Theorem \cite{EE} implies that every connected
component of $\symp^m(M,\omega)$ is contractible.
\smallskip\\
(Floer homology)
To every $\phi\in\symp^m(M,\omega)$ symplectic Floer homology
theory assigns a $\Z_2$-graded vector space $\HF(\phi)$ over $\Z_2$, with an
additional multiplicative structure, called the quantum cap product,
$
H^*(M;\Z_2)\otimes\HF(\phi)\To\HF(\phi).
$
For $\phi=id_M$ the symplectic Floer homology $\HF(id_M)$ are  canonically isomorphic to ordinary homology  $H_*(M;\Z_2)$ and quantum cap product agrees with the ordinary cap product.
Each $\psi\in\diff^+(M)$ induces an isomorphism
$\HF(\phi)\cong\HF(\psi^{-1}\phi\psi)$ of $H^*(M;\Z_2)$-modules. 
\smallskip\\
(Invariance) 
If $\phi,\phi'\in\symp^m(M,\omega)$ are isotopic, then
$\HF(\phi)$ and $\HF(\phi')$ are naturally isomorphic as
$H^*(M;\Z_2)$-modules. 
This is proven in \cite[Page 7]{S}. Note that every Hamiltonian perturbation
of $\phi$ (see \cite{ds}) is also in $\symp^m(M,\omega)$.
\smallskip\\ 
Now let $g$ be a mapping class of $M$, i.e. an isotopy class of $\diff^+(M)$. 
Pick an area form $\omega$ and a
representative $\phi\in\symp^m(M,\omega)$ of $g$. 
Then $\HF(\phi)$ is an invariant of $g$, which is
denoted by $\HF(g)$. Note that $\HF(g)$ is independent of the choice of an area
form $\omega$ by Moser's isotopy theorem \cite{Mo} and naturality of Floer homology.

\section{ Symplectic Floer homology\label{sec:floer}}

Let $\phi\in\symp(M,\omega)$.There are two ways of constructing Floer
homology detecting its fixed points, $\Fix(\phi)$. Firstly, the graph of $\phi$
is a Lagrangian submanifold of $M\times M,(-\omega)\times\omega)$ and its fixed points correspond to the intersection points of graph($\phi$) with the
diagonal $\Delta=\{(x,x)\in M\times M\}$. Thus we have the Floer homology of the Lagrangian intersection $\HF(M\times M,\Delta, graph (\phi))$.
This intersection is transversal if the fixed points of $\phi$ are nondegenerate, i.e. if 1 is not an eigenvalue of $d\phi(x)$, for $x\in\Fix(\phi)$.
The second approach was mentioned by Floer in \cite{Floer1} and presented with
details by Dostoglou and Salomon  in \cite{ds}.We follow here  Seidel's approach \cite{S} which, comparable with \cite {ds}, uses a larger
class of perturbations, but such that the perturbed action form is still
cohomologous to the unperturbed. As a consequence, the usual invariance of
Floer homology under Hamiltonian isotopies is extended to the stronger
property stated above.
Let now  $\phi\in\symp^m(M,\omega)$, i.e  $\phi$ is  monotone.  
Firstly, we give  the definition of $\HF(\phi)$  in the special case
where all the fixed points of $\phi$ are non-degenerate, i.e. for all $y\in\fix(\phi)$, $\det(\id-\de\phi_y)\ne0$, and then
following Seidel´s approach  \cite{S} we consider general case when
$\phi$ has degenerate fixed points.
 Let $\Of = \{ y \in C^{\infty}(\R,M)\,|\, y(t) = \phi(y(t+1)) \}$ be the twisted free loop space, which is also
the space of sections of $T_\phi \rightarrow S^1$. The action form is the
closed one-form $\alpha_\phi$ on $\Of$ defined by
$$
\alpha_\phi(y) Y = \int_0^1 \omega(dy/dt,Y(t))\,dt.
$$
where $y\in\Of$ and $Y\in T_y\Of$, i.e. $Y(t)\in T_{y(t)}M$ and
$Y(t)=\de\phi_{y(t+1)}Y(t+1)$ for all $t\in\R$.
 
The tangent bundle of any symplectic manifold admits an almost complex structure $ J:TM\To TM$ which is compatible with $\omega$ in sense that $(v,w)=\omega(v,Jw)$ defines a Riemannian metric.
 Let $J=(J_t)_{t \in \R}$ 
be a smooth path of $\omega$-compatible almost  complex structures on 
$M$ such that $J_{t+1}=\phi^*J_t$.
If $Y,Y'\in T_y\Of$, then 
$\int_0^1\omega(Y'(t),J_t Y(t))\de t$ defines a metric 
on the loop space $\Of$. So  the critical points of $\alpha_\omega$ are the constant paths in  $\Of$ and  hence the fixed points of $\phi$. The negative gradient lines of $\alpha_\omega$ with respect to the
metric above  are solutions of the partial differential equations
with boundary conditions
\begin{equation}\label{eq:corbit}
\left\{\begin{array}{l}
u(s,t) = \phi(u(s,t+1)), \\
\p_s u + J_t(u)\p_t u = 0, \\
\lim_{s\to\pm\infty}u(s,t) \in \Fix(\phi)
\end{array}\right.
\end{equation}
These are exactly  Gromov's pseudoholomorphic  curves \cite{Gromov}. 

For $y^\pm\in\fix(\phi)$, let $\M(y^-,y^+;J,\phi)$ denote the space of smooth maps $u:\R^2\to M$ which satisfy the  equations \eqref{eq:corbit}.
Now to every $u\in\M(y^-,y^+;J,\phi)$ is associated a Fredholm operator
$\De_u$ which linearizes (\ref{eq:corbit}) in suitable Sobolev spaces. The
index of this operator is given by the so called Maslov index $\mu(u)$,
which satisfies $\mu(u)=\deg(y^+)-\deg(y^-)\text{ mod }2$, where  $(-1)^{\deg y}=\sign(\det(\id-\de\phi_y))$. We have no bubling, since for surface 
$\pi_2(M)=0$. For a generic
$J$, every $u\in\M(y^-,y^+;J,\phi)$ is regular, meaning that $\De_u$ is onto. 
Hence, by the implicit function theorem, $\M_k(y^-,y^+;J,\phi)$ is
a smooth $k$-dimensional manifold, where $\M_k(y^-,y^+;J,\phi)$ denotes the
subset of those $u\in\M(y^-,y^+;J,\phi)$ with $\mu(u)=k\in\Z$.
Translation of the $s$-variable defines a free $\R$-action on 1-dimensional
manifold $\M_1(y^-,y^+;J,\phi)$ and hence the quotient is a discrete set of points. The energy of a map $u:\R^2\to M$ is given by $
E(u) = \int_{\R}\int_0^1 \omega\big(\p_tu(s,t),J_t\p_tu(s,t)\big)\,\de t\de s$  for all $y\in\fix(\phi)$.  P.Seidel has proved  in \cite{S} that if $\phi$ is monotone, then the energy is constant on each $\M_k(y^-,y^+;J,\phi)$.
Since all fixed points of $\phi$ are nondegenerate the set  $\fix(\phi)$ is a finite set and the
$\Z_2$-vector space  $
\CF(\phi) := \Z_2^{\#\fix(\phi)}$
admits a $\Z_2$-grading with $(-1)^{\deg y}=\sign(\det(\id-\de\phi_y))$, 
for all $y\in \fix(\phi)$.
The boundness of the energy  $E(u)$   for monotone  $\phi$  implies that   the  0-dimensional  quotients   $\M_1(y_-,y_+,J,\phi)/\R$   are actually finite sets. Denoting by  $n(y_-,y_+)$  the number of
points mod 2 in each of them, one defines a differential  $\partial_{J}:
CF_*(\phi) \rightarrow CF_{* + 1}(\phi)$  by $\partial_{J}y_- =
\sum_{y_+} n(y_-,y_+) {y_+}$.  Due to gluing theorem  this Floer boundary operator satisfies  $\partial_{J} \circ
\partial_{J} = 0$.  For gluing  theorem to hold one needs again the  boundness of the energy $E(u)$ .  It follows that  $ (\CF(\phi),\partial_{J})$  is a chain complex  and its homology is by definition the Floer homology of  $\phi$   denoted $HF_*(\phi)$. It  is independent of $J$ and is an invariant of $\phi$.

If $\phi$ has degenerate fixed points one needs to perturb equations 
\eqref{eq:corbit} in order to define the Floer homology. Equivalently, one
could say that the action form needs to be perturbed.
 The necessary analysis is given in \cite{S}  is essentially  the same as in the slightly different
situations considered in \cite{ds}. But  Seidel's approach also differs from the usual one in \cite{ds}. He uses a larger
class of perturbations, but such that the perturbed action form is still
cohomologous to the unperturbed.

\section{Nielsen numbers and Floer homology\label{sec:numbers}}

Before discussing the   results of the paper, we briefly describe the few
basic notions of Nielsen fixed point theory which will be used.
We assume  $X$ to be a connected, compact
polyhedron and $f:X\rightarrow X$ to be a continuous map.
Let $p:\tilde{X}\rightarrow X$ be the universal cover of $X$
and $\tilde{f}:\tilde{X}\rightarrow \tilde{X}$ a lifting
of $f$, i.e. $p\circ\tilde{f}=f\circ p$.
Two liftings $\tilde{f}$ and $\tilde{f}^\prime$ are called
{\sl conjugate} if there is a $\gamma\in\Gamma\cong\pi_1(X)$
such that $\tilde{f}^\prime = \gamma\circ\tilde{f}\circ\gamma^{-1}$.
The subset $p(Fix(\tilde{f}))\subset Fix(f)$ is called
{\sl the fixed point class of $f$ determined by the lifting class $[\tilde{f}]$}.Two fixed points $x_0$ and $x_1$ of $f$ belong to the same fixed point class iff  there is a path $c$ from $x_0$ to $x_1$ such that $c \cong f\circ c $ (homotopy relative endpoints). This fact can be considered as an equivalent definition of a non-empty fixed point class.
 Every map $f$  has only finitely many non-empty fixed point classes, each a compact  subset of $X$.
A fixed point class is called {\sl essential} if its index is nonzero.
The number of essential fixed point classes is called the {\sl Nielsen number}
of $f$, denoted by $N(f)$.The Nielsen number is always finite.
$R(f)$ and $N(f)$ are homotopy invariants.
In the category of compact, connected polyhedra, the Nielsen number
of a map is, apart from certain exceptional cases,
 equal to the least number of fixed points
 of maps with the same homotopy type as $f$.

\subsection{Periodic diffeomorphisms}

The following Lemma was  first proven in \cite{J}. We repeat here the
arguments from \cite{G}.

\begin{lemma}\label{lemma:jiang}\cite{J}
Let $\phi$  a non-trivial orientation
preserving periodic diffeomorphism of a compact connected surface $M$
of Euler characteristic $\chi(M) \leq0$. Then each fixed point class of $\phi$
consists of a single point. 
\end{lemma}

\begin{proof} 
First assume that $M$ is closed. The uniformization theorem states 
that in every conformal class of metrics on $M$, there is a unique metric of
constant curvature $-1$ if $\chi(M)<0$ or 0 if $\chi(M)=0$. This implies that
the unique representative of a $\phi$-invariant conformal class  of metrics
( such a class exists since $\phi$ is finite order) is itself $\phi$-invariant. Hence we can pick a $\phi$-invariant metric of constant curvature $-1$ or 0 on $M$ and 
lift $\phi$ to an isometry $\tilde{\phi}$ of the universal cover
$\tilde{M}$ of $M$. $\tilde{M}$ is either isometric to the hyperbolic plane
$H^2$ or the Euclidean plane $\R^2$. 

Let $x\in\fix(\phi)$ and 
let $\tilde{\phi},\tilde{x}$ be lifts of $\phi,x$ to
$\tilde{M}$, such that $\tilde{\phi}(\tilde{x})=\tilde{x}$. 
Note, that a fixed point of $\phi$ is in the same class as $x$ if and only
if it can be lifted 
to a fixed point of $\tilde{\phi}$. Assume by contradiction that
$\tilde{y}\ne\tilde{x}$ is a fixed point of $\tilde{\phi}$. It follows that
the unique geodesic going through $\tilde{x}$ and $\tilde{y}$ is pointwise fixed by
$\tilde{\phi}$. In particular, since $\tilde{\phi}$ preserves
orientation, $\de\tilde{\phi}_{\tilde{x}}=\id$. This implies
that $\tilde{\phi}=\id$, because an isometry of $H^2$ or $\R^2$ is
determined by its value and differential at one point. This proves claim 3 in
the case that $M$ is closed.  

The case $\p M\ne\emptyset$ is reduced to the above case by gluing two
copies of $M$ together along a $\phi$-invariant tubular neighborhood of $\p
M$. The glued manifold is closed and of Euler characteristic $\leq0$;
$\varphi$ extends to a non-trivial diffeomorphism $\phi'$, which is
orientation preserving and of finite order.
Hence, every fixed point class of $\phi'$ is a single point.
The same therefore holds for $\phi$. This ends the proof. 
\end{proof}

R. Gautschi gave two criteria for monotonicity which we use later on.

 Let $\omega$ be an area form on $M$ and
$\phi\in\symp(\Sigma,\omega)$. 
\begin{lemma}\label{lemma:monotone1}\cite{G}
Assume that every class $\alpha\in\ker(\id-\phi_*)\subset H_{1}(M;\Z)$ is
represented by a map $\gamma:S\to\fix(\phi)$, where $S$ is a compact oriented
1-manifold. Then $\phi$ is monotone.
\end{lemma}
\begin{lemma}\label{lemma:monotone2}\cite{G}
If $\phi^k$ is monotone for some $k>0$, then $\phi$ is monotone. 
If $\phi$ is monotone, then $\phi^k$ is monotone for all $k>0$. 
\end{lemma}
\begin{proof} We repeat  Gautschi arguments  from \cite{G} here.
Recall that $\tf$ is the orbit space of the $\Z$-action
$n\cdot(t,x) = (t+n,\phi^{-n}(x))$, where $n\in\Z$ and $(t,x)\in\R\times\Sigma$. If we only divide out by the 
subgroup $k\Z$, for $k\in\N_{>0}$, we naturally get the mapping torus of
$\phi^k$. 
Further dividing by $\Z/k\Z$ defines the $k$-fold covering map
$p_k:T_{\phi^k}\to\tf$. 
It is straight forward to check that 
\begin{equation}
p_k^*[\omega_\phi] = [\omega_{\phi^k}] \quad\text{and}\quad 
p_k^*c_\phi = c_{\phi^k}.
\label{eq:iterate}\end{equation}
The first equality follows immediately from the definitions. 
To prove the second, note that
$  
p_k^*\big((TM\times\R)/\Z\big) \cong (TM\times\R)/k\Z \cong V_{\phi^k}, 
$
where the $\Z$-action on $\R\times T\Sigma$ is given by 
$n\cdot(t,\xi_x)=(t+n,\de\phi_x^{-n}\xi_x)$, for $n\in\Z$ and $\xi_x\in
T_x M$. 
The lemma follows from \eqref{eq:iterate} and the fact that $p_k^*$ is
injective. To prove injectivity, define the map
$a^k_*:H^2(T_{\phi^k};\R)\to H^2(T_\phi;\R)$ by averaging differential forms;
$a^k_*$ is a left inverse of $p_k^*$, i.e. $a^k_*\circ p_k^*=\id$. This ends
the proof of the lemma.
\end{proof}
We shall say  that $\phi:M\rightarrow M$ is a  periodic map of  period $m$,
if $\phi^m$ is  the identity map $id_M:M\rightarrow M$.

\begin{theorem}\label{thm:main}
If $\phi$ is  a non-trivial orientation
preserving periodic diffeomorphism of a compact connected surface $M$
of Euler characteristic $\chi(M) <  0$, then $\phi$ is monotone with
respect to some $\phi$-invariant area form and 
$$
\HF(\phi) \cong
\Z_2^{N (\phi)},  \,\,\,  \dim\HF(\phi)= N (\phi),
$$
where  $ N (\phi)$ denotes the Nielsen  number of  $\phi$.
\end{theorem}
\begin{proof}
Let  $\phi$ be  a periodic diffeomorphism  of least period $l$.
 First note that
if $\tilde{\omega}$ is an area form on $M$, then area form
$\omega:=\sum_{i=1}^\ell(\phi^i)^*\tilde{\omega}$ is 
 $\phi$-invariant, i.e. $\phi\in\symp(M,\omega)$. By periodicity of  $\phi$,
 $\phi^l$ is  the identity map $id_M:M\rightarrow M$.
Then from Lemmas \ref{lemma:monotone1}  and \ref{lemma:monotone2} it follows that $\omega$ can be chosen such that $\phi\in\symp^m(M,\omega)$.

 Lemma \ref{lemma:jiang}  implies that every
$y\in\fix(\phi)$ forms a different fixed point class of
$\phi$, so $ \#\fix(\phi)= N (\phi)$.
This has an immediate consequence for the Floer complex $(\CF(\phi),\p_{J})$
with respect to a generic $J=(J_t)_{t\in\R}$.
If $y^\pm\in\fix(\phi)$ are in different fixed point classes, then 
$\M(y^-,y^+;J,\phi)=\emptyset$. This follows from the first equation in
\eqref{eq:corbit}.
Then the boundary map in  the Floer complex is zero $\partial_{J}=0$ and 
$\Z_2$-vector space  $
\CF(\phi) := \Z_2^{\#\fix(\phi)}=\Z_2^{N(\phi)}$. This immediately  implies
$
\HF(\phi) \cong
\Z_2^{N (\phi)}
$ and $ \dim\HF(\phi)= N (\phi)$.

\end{proof}

\subsection{Algebraically finite mapping classes}

 A mapping class of $M$ is called  algebraically finite if it
does not have any pseudo-Anosov components in the sense of Thurston's
theory of surface diffeomorphism.The term algebraically finite goes back to J. Nielsen.\\
 R. Gautschi \cite{G} defined the notion of a
diffeomorphism of finite type  for surface diffeomorphisms. These are
special representatives of algebraically finite mapping classes
adopted to the symplectic geometry. 
\begin{dfn}\label{def:ftype}\cite{G}
We call $\phi\in\diff_+(M)$ of {\bf finite type} if the following holds.
There is a $\phi$-invariant finite union $N\subset M$ of disjoint
non-contractible annuli such that:
\smallskip\\
(1) $\phi|M \setminus N$ is periodic, i.e. there exists
$\ell>0$ such that $\phi^\ell|M \setminus N=\id$.
\smallskip\\
(2) Let $N'$ be a connected component of $N$ and $\ell'>0$ be the
smallest integer such that $\phi^{\ell'}$ maps $N'$ to itself. Then
$\phi^{\ell'}|N'$ is given by one of the following two models with respect to
some coordinates $(q,p)\in I\times S^1$:
\medskip\\
\begin{minipage}{4cm}
(twist map)
\end{minipage}
\begin{minipage}{6cm}
$(q,p)\Mapsto(q,p-f(q))$
\end{minipage}
\medskip\\
\begin{minipage}{4cm}
(flip-twist map)
\end{minipage}
\begin{minipage}{6cm}
$(q,p)\Mapsto(1-q,-p-f(q))$,
\end{minipage}
\medskip\\
where $f:I\to \R$ is smooth and strictly monotone. 
A twist map is called  positive or negative,
if $f$ is increasing or decreasing. 
\smallskip\\
(3) Let $N'$ and $\ell'$ be as in (2).
If $\ell'=1$ and $\phi|N'$ is a twist map, then $\im(f)\subset[0,1]$,
i.e. $\phi|\text{int}(N')$ has no fixed points.
\smallskip\\
(4) If two connected components of $N$ are homotopic, then the corresponding local
models of $\phi$ are either both positive or both negative twists.

\end{dfn}
The term flip-twist map is taken from \cite{JG}.

By $M_{id}$ we denote the union of the components
of $M\setminus\text{int}(N)$, where $\phi$ restricts to the identity.

The next  lemma describes the set of fixed point classes of $\phi$. It
is a special case of a theorem by B.~Jiang and J.~Guo~\cite{JG}, which gives
for any mapping class a representative that realizes its
Nielsen number.
\begin{lemma}[Fixed point classes]\label{lemma:fclass}\cite{JG}
Each fixed point class of $\phi$ is either a connected component of $M_{id}$ or
consists of a single fixed point. A fixed point $x$ of the second type satisfies
$\det(\id-\de\phi_x)>0$.
\end{lemma}

The monotonicity of diffeomorphisms of finite type
was  investigated in details in the recent preprint  by R. Gautschi\cite{G}.
We describe now his results.
Let $\phi$ be a diffeomorphism of finite type and
$\ell$ be as in (1).
Then $\phi^\ell$ is the product of (multiple)  Dehn twists along $N$.
Moreover, two parallel Dehn twists have the same sign, by (4). We say that
$\phi$ has  uniform twists, if $\phi^\ell$ is the product of only
positive, or only negative Dehn twists.
\smallskip\\
Furthermore, we denote by $\ell$ the smallest positive integer such that
$\phi^\ell$ restricts to the identity on $M\setminus N$.

If $\omega'$ is an area form on $M$ which is the standard form
$\de q\wedge\de p$ with respect to the $(q,p)$-coordinates on $N$, then
$\omega:=\sum_{i=1}^\ell(\phi^i)^*\omega'$ is
standard on $N$ and $\phi$-invariant, i.e. $\phi\in\symp(M,\omega)$.
To prove that $\omega$ can be chosen such that $\phi\in\symp^m(M,\omega)$,
Gautschi distinguishes two cases: uniform and non-uniform twists. In the first case he proves the following stronger statement.
\begin{lemma}\label{lemma:monotone3}\cite{G}
If $\phi$ has uniform twists and $\omega$ is a $\phi$-invariant area
form, then $\phi\in\symp^m(M,\omega)$.
\end{lemma}

In the non-uniform case, monotonicity  does not
hold for arbitrary $\phi$-invariant area forms.
\begin{lemma}\label{lemma:monotone4}\cite{G}
If $\phi$ does not have uniform twists, there exists a $\phi$-invariant
area form $\omega$ such that $\phi\in\symp^m(M,\omega)$. Moreover,
$\omega$ can be chosen such that it is the standard form $\de q\wedge\de p$ on
$N$.
\end{lemma}

\begin{theorem}\label{thm:main0}
If $\phi$ is  a  diffeomorphism of finite type  of a compact connected surface $M$
of Euler characteristic $\chi(M) <  0$ and if  $\phi$ has  only isolated fixed points , then $\phi$ is monotone with
respect to some $\phi$-invariant area form and 
$$
\HF(\phi) \cong
\Z_2^{N (\phi)},  \,\,\,  \dim\HF(\phi)= N (\phi),
$$
where  $ N (\phi)$ denotes the Nielsen  number of  $\phi$.
\end{theorem}
\begin{proof}
From Lemmas \ref{lemma:monotone3}  and \ref{lemma:monotone4} it follows that $\omega$ can be chosen such that $\phi\in\symp^m(M,\omega)$.
 Lemma \ref{lemma:fclass}  implies that every
$y\in\fix(\phi)$ forms a different fixed point class of
$\phi$, so $ \#\fix(\phi)= N (\phi)$.
This has an immediate consequence for the Floer complex $(\CF(\phi),\p_{J})$
with respect to a generic $J=(J_t)_{t\in\R}$.
If $y^\pm\in\fix(\phi)$ are in different fixed point classes, then 
$\M(y^-,y^+;J,\phi)=\emptyset$. This follows from the first equation in
\eqref{eq:corbit}.
Then the boundary map in  the Floer complex is zero $\partial_{J}=0$ and 
$\Z_2$-vector space  $
\CF(\phi) := \Z_2^{\#\fix(\phi)}=\Z_2^{N(\phi)}$. This immediately  implies
$
\HF(\phi) \cong
\Z_2^{N (\phi)}
$ and $ \dim\HF(\phi)= N (\phi)$.

\end{proof}

\begin{rk}\label{rk:G}
R. Gautschi  has proved in preprint  \cite{G} that,
if $\phi$ is a diffeomorphism of finite type, then $\phi$ is monotone with
respect to some $\phi$-invariant area form and
$$
\HF(\phi) =
H_*(M_{id},\p_{M_{id}};\Z_2)\oplus
\Z_2^{L(\phi|M\setminus M_{id})}.
$$
Here, $L$ denotes the Lefschetz number( see section \ref{sec:zeta}).

In  the  theorem \ref{thm:main0} the  set $ M_{id}$ is empty and  every fixed point of $\phi$ has fixed point index 1 \cite{JG}. 
The Lefschetz fixed point formula implies that 
$\#\fix(\phi)=N(\phi)=L(\phi)$ .
So, theorem \ref{thm:main0} now  follows also  from   result of R.Gautschi.
\end{rk}

\subsection{Hyperbolic diffeomorphisms of 2-dimensional  torus}

\begin{theorem}\label{thm:main1}
If $\phi$ is a  hyperbolic  diffeomorphism of a 2-dimensional torus $T^2$,
 then $\phi$ is symplectic  and 
$$
\HF(\phi) \cong
\Z_2^{N (\phi)},  \,\,\,  \dim\HF(\phi)= N (\phi)
$$
where  $ N (\phi)=|\det (E-\phi_*)|$ denotes the Nielsen  number of  $\phi$
and $\phi_*$ is an  induced homomorphism on the fundamental group of
    $T^2$.
\end{theorem}
\begin{proof}
Hyperbolicity of $\phi$  means that the covering linear map
 $\tilde \phi  : R^2 \rightarrow R^2 $, $\det\tilde \phi=1$  has no eigenvalue of modulus one. The hyperbolic  diffeomorphism of a 2-dimensional torus $T^2$ is area preserving so symplectic.
In fact, the covering map $\tilde{\phi} $ has a unique fixed point, which is the origin;
hence, by the covering homotopy theorem , the fixed points of $\phi$ are
pairwise Nielsen  nonequivalent.The index of each Nielsen  fixed point  class, consisting of one fixed
point, coincides with its Lefschetz index, and by the hyperbolicity of $\phi$, the later
is not equal to zero.Thus the Nielsen number $ N(\phi)= \#\fix(\phi)$.
If $y^\pm\in\fix(\phi)$ are in different Nielsen  fixed point classes, then 
$\M(y^-,y^+;J,\phi)=\emptyset$. This follows from the first equation in
\eqref{eq:corbit}.
Then the boundary map in  the Floer complex is zero $\partial_{J}=0$ and 
$\Z_2$-vector space  $
\CF(\phi) := \Z_2^{\#\fix(\phi)}=\Z_2^{N(\phi)}$. This immediately  implies
$
\HF(\phi) \cong
\Z_2^{N (\phi)}
$ and $ \dim\HF(\phi)= N (\phi)$.

\end{proof}
\begin{rk}
It is interesting to compare this result with
the first computation by Marcin Po\'zniak 
of the Floer homology of linear symplectomorphisms in case 
of  torus\cite{Po}.

\end{rk}

\section{Symplectic zeta functions\label{sec:zeta}}

Let $\Gamma = \pi_0(Diff^+(M))$ be the mapping class group of a closed connected oriented surface $M$ of genus $\geq 2$. Pick an everywhere positive two-form $\omega$ on $M$. A isotopy  theorem of Moser \cite{Mo} says that each  mapping class of  $g \in \Gamma$,  i.e. an isotopy class of $Diff^+(M)$,  admits representatives which preserve $\omega$. Due  to Seidel\cite{S} we can
pick  a monotone  representative $\phi\in\symp^m(M,\omega)$ of $g$. 
Then $\HF(\phi)$ is an invariant of $g$, which is
denoted by $\HF(g)$. Note that $\HF(g)$ is independent of the choice of an area
form $\omega$ by Moser's  theorem  and naturality of Floer homology. 
By Gautschi lemma  \ref{lemma:monotone2}  symplectomorphisms  $\phi^n$
are  also monotone for all $n>0$.
Taking a dynamical point of view,
 we consider the iterates of monotone symplectomorphism  $\phi$
and  define the first symplectic  zeta function of $\phi$
 as the following power series:

\begin{eqnarray*}
 \chi_\phi(z) = 
 \exp\left(\sum_{n=1}^\infty \frac{\chi(\HF(\phi^n))}{n} z^n \right),
\end{eqnarray*}
where $\chi(\HF(\phi^n))$ is the Euler characteristic of Floer homology group
of $\phi^n$. 
Then $ \chi_\phi(z)$ is an invariant of $g$, which we
denote by $ \chi_g(z)$.

Let us consider  the  Lefschetz  zeta function
 $$
  L_\phi(z) := \exp\left(\sum_{n=1}^\infty \frac{L(\phi^n)}{n} z^n \right),
$$
   where
 $
   L(\phi^n) := \sum_{k=0}^{\dim X} (-1)^k \tr\Big[\phi_{*k}^n:H_k(M;\Q)\to H_k(M;\Q)\Big]
 $
 is the Lefschetz number of $\phi^n$. The Lefschetz zeta function is always a rational function of $z$ and
is given by the formula:
$$
 L_\phi(z) = \prod_{k=0}^{\dim X}
          \det\big(E-\phi_{*k}.z\big)^{(-1)^{k+1}}.
$$

\begin{theorem}\label{thm:lef}
Symplectic zeta function $ \chi_\phi(z)$ is a rational function of $z$ and 
$$
 \chi_\phi(z)= L_\phi(z)=\prod_{k=0}^{\dim X}
          \det\big(E-\phi_{*k}.z\big)^{(-1)^{k+1}}.
$$
\end{theorem}

\begin{proof} 
If for every $n$ all the fixed points of $\phi^n$ are non-degenerate, i.e. for all $x\in\fix(\phi^n)$, $\det(\id-\de\phi^n(x))\ne0$, then we have( see section 
\ref{sec:floer}): 
$$\chi(\HF(\phi^n))=\sum_{x=\phi^n(x)} \sign(\det(\id-\de\phi^n(x)))=L(\phi^n).
$$
If we have degenerate fixed points one needs   to perturb equations 
\eqref{eq:corbit} in order to define the Floer homology.
 The necessary analysis is given in \cite{S}  is essentially  the same as in the slightly different situations considered in \cite{ds}, where the  above connection between the Euler characteristic and the Lefschetz number 
 was firstly  established.
\end{proof}

\begin{rk} Theorem \ref{thm:lef} shows that symplectic 
zeta function $ \chi_\phi(z)$  counts symplectic periodic points of $\phi$ algebraically-
in the  Lefschetz way.
\end{rk}
 
The next issue is a relation of the symplectic zeta function with
the Reidemeister torsion, a very important topological invariant. 
 We will show that special value of symplectic zeta function
$ \chi_\phi(z)$ is a Reidemeister torsion.
 Let $T_\phi := (X\times I)/(x,0)\sim(\phi(x),1)$ be the
mapping tori of  $\phi$.
We shall consider the bundle $p:T_\phi\rightarrow S^1$
over the circle $S^1$ with fibers $M$.
We assume here that $E$ is a flat, complex vector bundle with
finite dimensional fibre and base $S^1$. We form its pullback $p^*E$
over $T_\phi$.
Note that the vector spaces $H^i(p^{-1}(b),\C)$ with
$b\in S^1$ form a flat vector bundle over $S^1$,
which we denote $H^i M$. The integral lattice in
$H^i(p^{-1}(b),\R)$ determines a flat density by the condition
that the covolume of the lattice is $1$.
We suppose that the bundle $E\otimes H^i M$ is acyclic for all
$i$. Under these conditions D. Fried \cite{fri1} has shown that the bundle
$p^* E$ is acyclic, and we have
\begin{equation}\label{eq:tor}
 \tau(T_\phi;p^* E) = \prod_i \tau(S^1;E\otimes H^i M)^{(-1)^i}.
\end{equation}
Let $g$ be the prefered generator of the group
$\pi_1 (S^1)$ and let $A=\rho(g)$ where
$\rho:\pi_1 (S^1)\rightarrow GL(V)$.
Then the holonomy around $g$ of the bundle $E\otimes H^i M$
is $A\otimes \phi^*_i$.
 
Since $\tau(S^1;E)=\mid\det(I-A)\mid$ it follows from \eqref{eq:tor}
that
\begin{equation}
 \tau(T_\phi;p^* E) = \prod_i \mid\det(I-A\otimes \phi^*_i)\mid^{(-1)^i}.
\end{equation}
We now consider the special case in which $E$ is one-dimensional,
so $A$ is just a complex scalar $\lambda$ of modulus one.
Then in terms of the rational function $L_\phi(z)$ we have \cite{fri1}:
\begin{equation}\label{eq:tor1}
 \tau(T_\phi;p^* E) = \prod_i \mid\det(I-\lambda \cdot \phi^*_i)\mid^{(-1)^i}
             = \mid L_\phi(\lambda)\mid^{-1}
\end{equation}
This proves the following
\begin{theorem} The Reidemeister torsion is the special value of
the symplectic zeta function:
$$
 \tau\left(T_{\phi};p^*E\right)
 =
 \mid \chi_{\phi}(\lambda) \mid^{-1},
 $$
where $\lambda$ is the holonomy of the one-dimensional
flat complex bundle $E$ over $S^1$.
\end{theorem}

Now we define  the second symplectic zeta function for 
 monotone symplectomorphism  $\phi$  as the following power series:

\begin{eqnarray*}
 F_\phi(z) &= &
 \exp\left(\sum_{n=1}^\infty \frac{\dim\HF(\phi^n)}{n} z^n \right).
\end{eqnarray*}
Then $ F_\phi(z)$ is an invariant of mapping class  $g$, which we
denote by $ F_g(z)$.

Motivation for this definition is the  theorem \ref{thm:main} and nice analytical
properties of the Nielsen zeta function $N_\phi(z) = 
 \exp\left(\sum_{n=1}^\infty \frac{N(\phi^n)}{n} z^n \right)$, see \cite{f1,f2}. 
We denote the  numbers  $\dim\HF(\phi^n) $ by $N_n$. Let $ \mu(d), d \in N$,
be the M\"obius function from number theory. As is known, it is given
by the following equations:  $\mu(d)=0 $ if $d$ is divisible by a square different from one ;
$\mu(d)=(-1)^k $ if
 $d$ is not divisible  by a square different from one , where $k$ denotes the number of
 prime divisors of $d$; $ \mu(1)=1$.

 \begin{theorem}
Let $\phi$ be  a non-trivial orientation
preserving periodic diffeomorphism  of least period $m$ of a compact connected surface $M$
of Euler characteristic $\chi(M) < 0$  . Then the 
  zeta function $ F_\phi(z)$ is  a radical of a rational function and 
$$
 F_\phi(z) =\prod_{d\mid m}\sqrt[d]{(1-z^d)^{-P(d)}},
$$
 where the product is taken over all divisors $d$ of the period $m$, and $P(d)$ is the integer
$$  P(d) = \sum_{d_1\mid d} \mu(d_1)N_{d\mid d_1} .  $$
\end{theorem}
\begin{proof}
Since $\phi^m = id $,then  for each $j, N_j=N_{m+j}$. If  $(k,m)=1$, then  there exist positive integers $t$
and $q$ such that $kt=mq+1$. So $ (\phi^k)^t=\phi^{kt}= \phi^{mq+1}=\phi^{mq}\phi=(\phi^m)^{q}\phi =  \phi$.
Consequently, $  \Fix((\phi^k)^t)=\Fix(\phi) $. 
We have $\Fix(\phi)\subset \Fix((\phi^k)$ and $\Fix((\phi^k)\subset \Fix((\phi^k)^t)=\Fix(\phi) $. Then  $\Fix(\phi) = \Fix((\phi^k)$ and 
$ N_1=N_k $. One can prove completely analogously that $ N_d= N_{di} $, if $(i, m/d) =1$,
where $d$ is a divisor of $m$. Using these sequences of equal  numbers, one can regroup
the terms of the series in the exponential of the  zeta function so as to get logarithmic
functions by adding and subtracting missing terms with necessary coefficients.
We show how to do this first for period $m=p^l$, where $p$ is a prime number.  We have the
following sequence of equal  numbers:
 $$
 N_1=N_k, (k,p^l)=1 \mbox{(i.e., no}\, N_{ip}, N_{ip^2},\dots
, N_{ip^l},  i=1,2,3,\dots),
 $$
 $$ 
 N_p=N_{2p}=N_{3p}=\dots =N_{(p-1)p}=N_{(p+1)p}= \dots \mbox{ (no} \,
N_{ip^2}, N_{ip^3}, \dots  , N_{ip^l} ) 
$$
etc.; finally,
$$
N_{p^{l-1}}=N_{2p^{l-1}}=\dots \mbox{(no }\, N_{ip^l})
$$
and separately the number $N_{p^l}$.\\
 Further,
\begin{eqnarray*}
\sum_{i=1}^\infty \frac{N_i}{i} z^i & = & \sum_{i=1}^\infty \frac{N_1}{i} z^i
       +\sum_{i=1}^\infty \frac{(N_p -N_1)}{p}\frac{ {z^p}^i}{i} + \\
                                                    & +  &\sum_{i=1}^\infty \frac{(N_{p^2} -(N_p -N_1)- N_1)}{p^2} \frac{{z^{p^2}}^i}{i} + ...\\
 & + & \sum_{i=1}^\infty \frac{(N_{p^l} - ...- (N_p -N_1)- N_1)}{p^l}\frac{ {z^{p^l}}^i}{i}\\
 &=& -N_1\cdot \log (1-z) + \frac{N_1-N_p}{p}\cdot
                                                    \log (1-z^p) +\\
& + & \frac{N_p-N_{p^2}}{p^2}\cdot \log (1-z^{p^2}) + ...\\
& + & \frac{N_{p^{l-1}}-N_{p^l}}{p^l}\cdot \log (1-z^{p^l}).
\end{eqnarray*}
For an arbitrary period $m$ , we get completely analogously,

\begin{eqnarray*}
F_f(z) & = & \exp\left(\sum_{i=1}^\infty \frac{N_i}{i} z^i \right)\\
           & = & \exp\left(\sum_{d\mid m} \sum _{i=1}^\infty \frac{P(d)}{d}\cdot\frac{{z^d}^i}{i}\right)\\
           & = & \exp\left(\sum_{d\mid m}\frac{P(d)}{d}\cdot \log (1-z^d)\right)\\
           & = & \prod_{d\mid m}\sqrt[d]{(1-z^d)^{-P(d)}},
\end{eqnarray*}
where the integers $P(d)$ are calculated recursively by the formula
$$
P(d)= N_d  - \sum_{d_1\mid d; d_1\not=d} P(d_1).
$$
Moreover, if the last formula is rewritten in the form
$$
N_d=\sum_{d_1\mid d}\mu(d_1)\cdot P(d_1)
$$
and one uses  the M\"obius Inversion law for real functions in number theory, then
$$
P(d)=\sum_{d_1\mid d}\mu(d_1)\cdot N_{d/d_1},
$$
where $\mu(d_1)$ is the M\"obius function in number theory. The lemma is proved.
\end{proof}
\begin{cor}
If in Theorem 2 the period $m$ is a prime number, then
$$
F_f(z)  =   \frac{1}{(1-z)^{N_1}}\cdot \sqrt[m]{(1-z^m)^{N_1 - N_m}}.
$$
\end{cor}
We denote by $\zeta_\phi(z)$ the  Artin-Mazur zeta function 
$$
 \zeta_\phi(z)  :=  \exp\left(\sum_{n=1}^\infty \frac{I(\phi^n)}{n} z^n \right),
$$ where $I(\phi^n)$ is the number of isolated fixed points of $\phi^n$.

\begin{theorem}\label{thm:main2}
If  $\phi$ is a  hyperbolic  diffeomorphism of a 2-dimensional torus $T^2$,
then the  symplectic zeta function $ F_\phi(z)$ is  a rational function
and $F_\phi(z)=N_\phi(z)=\zeta_\phi(z)=(L_\phi(\sigma\cdot z))^{(-1)^r}$ , where  $r$ is equal to the number of $\lambda_i \in Spec(\tilde \phi) $ such that $ \mid \lambda_i \mid > 1$, $p$ is equal to the number of $\mu_i \in Spec(\tilde \phi)$ such that $\mu_i <-1$ and  $\sigma=(-1)^p.$
\end{theorem}
\begin{proof}
From theorem \ref{thm:main1} and  \cite {bbpt}  it follows that $\dim\HF(\phi^n)=\#\fix(\phi^n)=I(\phi^n)= N (\phi^n)=\mid \det(E-\tilde {\phi^n}) \mid =\mid L(\phi^n) \mid $ \cite {bbpt}. Thus
 $\dim\HF(\phi^n)=\#\fix(\phi^n)=I(\phi^n)= N(\phi^n)=(-1)^{r+pn}\cdot \det(E-\tilde{\phi^n})=(-1)^{r+pn}\cdot L(\phi^n)$. Now a direct 
computation ends the proof of the theorem.

\end{proof}

\section{Asymptotic  invariant. Concluding  remarks and conjectures\label{sec:problem}}

\subsection{Topological entropy and the Nielsen numbers}

The most widely used measure for the complexity of a dynamical system is the topological
entropy. For the convenience of the reader, we include its definition.
 Let $ f: X \rightarrow X $ be a self-map of a compact metric space. For given $\epsilon > 0 $
 and $ n \in N $, a subset $E \subset X$ is said to be $(n,\epsilon)$-separated under $f$ if for
 each pair $x \not= y$ in $E$ there is $0 \leq i <n $ such that $ d(f^i(x), f^i(y)) > \epsilon$.
 Let $s_n(\epsilon,f)$  denote the largest cardinality of any $(n,\epsilon)$-separated subset $E$
 under $f$. Thus  $s_n(\epsilon,f)$ is the greatest number of orbit segments ${x,f(x),...,f^{n-1}(x)}$
 of length $n$ that can be distinguished one from another provided we can only distinguish
 between points of $X$ that are  at least $\epsilon$ apart. Now let
 $$
 h(f,\epsilon):= \limsup_{n} \frac{1}{n}\cdot\log \,s_n(\epsilon,f)
 $$
 $$
 h(f):=\limsup_{\epsilon \rightarrow 0} h(f,\epsilon).
 $$
 The number $0\leq h(f) \leq \infty $, which to be independent of the metric $d$ used, is called the topological entropy of $f$.
 If $ h(f,\epsilon)> 0$ then, up to resolution $ \epsilon >0$, the number $s_n(\epsilon,f)$ of
 distinguishable orbit segments of length $n$ grows exponentially with $n$. So $h(f)$
 measures the growth rate in $n$ of the number of orbit segments of length $n$
 with arbitrarily fine resolution.
   A basic relation between Nielsen numbers and topological entropy was found by N. Ivanov
   \cite{i1}. We present here a very short proof of Jiang  of the Ivanov's  inequality.
\begin{lemma}\label{lemma:ent}\cite{i1}
$$
h(f) \geq \limsup_{n} \frac{1}{n}\cdot\log N(f^n)
 $$
\end{lemma}
{\sc Proof}
 Let $\delta$ be such that every loop in $X$ of diameter $ < 2\delta $ is contractible.
 Let $\epsilon >0$ be a smaller number such that $d(f(x),f(y)) < \delta $ whenever $ d(x,y)<2\epsilon $. Let $E_n \subset X $ be a set consisting of one point from each essential fixed point class of $f^n$. Thus $ \mid E_n \mid =N(f^n) $. By the definition of $h(f)$, it suffices
 to show that $E_n$ is $(n,\epsilon)$-separated.
 Suppose it is not so. Then there would be two points $x\not=y \in E_n$ such that $ d(f^i(x), f^i(y)) \leq \epsilon$ for $o\leq i< n$ hence for all $i\geq 0$. Pick a path $c_i$ from $f^i(x)$ to
 $f^i(y)$ of diameter $< 2\epsilon$ for $ o\leq i< n$ and let $c_n=c_0$. By the choice of $\delta$
 and $\epsilon$ ,  $f\circ c_i \simeq c_{i+1} $ for all $i$, so $f^n\circ c_0\simeq c_n=c_0$. such that 
 This means $x,y$ in the same fixed point class of $f^n$, contradicting the construction of $E_n$.

 This inequality is remarkable in that it does not require smoothness of the map and provides a common lower bound for the topological entropy of all maps in a homotopy class.

We recall Thurston classification theorem for homeomorphisms of surfase $M$
of genus $\geq 2$.

\begin{theorem}\label{thm:thur}\cite{Th}
Every homeomorphism $\phi: M\rightarrow M $ is isotopic to a homeomorphism $f$
such that either\\
(1) $f$ is a periodic map; or\\
(2) $f$ is a pseudo-Anosov map, i.e. there is a number $\lambda >1$(stretching factor) and a pair of transverse measured foliations $(F^s,\mu^s)$ and $(F^u,\mu^u)$ such that $f(F^s,\mu^s)=(F^s,\frac{1}{\lambda}\mu^s)$ and $f(F^u,\mu^u)=(F^u,\lambda\mu^u)$; or\\
(3)$f$ is reducible map, i.e. there is a system of disjoint simple closed curves $\gamma=\{\gamma_1,......,\gamma_k\}$ in $int M$ such that $\gamma$ is invariant by $f$(but $\gamma_i$ may be permuted) and $ \gamma$ has a $f$-invariant 
tubular neighborhood $U$  such that each component of $M\setminus U$  has negative 
Euler characteristic and on each(not necessarily connected) $f$-component of
$M\setminus U$, $f$ satisfies (1) or (2).

\end{theorem}
The map $f$ above is called the Thurston canonical form of $f$. In (3) it 
can be chosen so that some iterate $f^m$ is a generalised Dehn twist on $U$.
Such a $f$ , as well as the $f$ in (1) or (2), will be called standard.
A key observation is that if $f$ is standard, so are all iterates of $f$.

\begin{lemma}\label{lem:flp}\cite{flp} Let $f$ be a pseudo-Anosov  homeomorphism with stretching factor $\lambda >1$ of surfase $M$
of genus $\geq 2$. Then
$$h(f)=log(\lambda)= \limsup_{n} \frac{1}{n}\cdot\log N(f^n)$$
\end{lemma}

\begin{lemma}\label{lem:j1}\cite{j1}
Suppose $f$ is a standard homeomorphism  of surfase $M$
of genus $\geq 2$ and  $\lambda$ is the largest stretching factor  of the  pseudo-Anosov pieces( $\lambda=1$ if there is no pseudo-Anosov piece).
Then
$$h(f)=log(\lambda)= \limsup_{n} \frac{1}{n}\cdot\log N(f^n)$$
\end{lemma}
\subsection{Asymptotic invariant}

The growth rate of a sequence $a_n$ of complex numbers is defined
by
 $$
 \grow ( a_n):= max \{1,  \limsup_{n \rightarrow \infty} |a_n|^{1/n}\}
$$ 

which could be infinity. Note that $\grow(a_n) \geq 1$ even if all $a_n =0$.
When  $\grow(a_n) > 1$, we say that the sequence $a_n$ grows exponentially.
\begin{dfn}
 We define the asymptotic invariant $ F^{\infty}(g)$ of  mapping class   $g \in \Gamma = \pi_0(Diff^+(M))$ to be the growth rate of the sequence
$\{a_n=\dim\HF(\phi^n)\}$ for  a monotone  representative $\phi\in\symp^m(M,\omega)$ of $g$:
$$ F^{\infty}(g):=\grow(\dim\HF(\phi^n)) $$
\end{dfn}
\begin{ex}
If $\phi$ is  a non-trivial orientation
preserving periodic diffeomorphism of a compact connected surface $M$
of Euler characteristic $\chi(M) <  0$ , then the  periodicity of the
sequence $\dim\HF(\phi^n)$ implies that  for the  corresponding  mapping class $g$ the  asymptotic invariant
$$ F^{\infty}(g):=\grow(\dim\HF(\phi^n))=1 $$

\end{ex}

\begin{ex} Let $\phi$ be  a monotone  diffeomorphism of finite type of a compact connected surface $M$
of Euler characteristic $\chi(M) <  0$ and $g$ a corresponding algebraically finite mapping class.
Let $U$ be the open regular neighborhood of the $k$ reducing
curves $\gamma_1,......,\gamma_k$ in the Thurston theorem, and $M_j$
be the component of $M\setminus U$.Let $F$ be a fixed point class of $\phi$.
Observe from \cite{JG} that if $F\subset M_j$, then $\ind(F, \phi)=\ind(F,\phi_j)$.So if $F$ is counted in $N(\phi)$ but not counted in $\sum_{j} N(\phi_j)$
, it must intersect $U$. But we see from \cite{JG} that a component
of $U$ can intersect at most 2 essential fixed point classes of $\phi$.
Hence we have $N(\phi)\leq \sum_{j} N(\phi_j)$.For the  monotone  diffeomorphism of finite type $\phi$ maps $\phi_j$ are  periodic. 
Applying last  inequality to $\phi^n$ and using  remark   \ref{rk:G} we have

$$
0\leq\dim\HF(\phi^n)= \dim H_*(M^{(n)}_{id},\p{M^{(n)}_{id}};\Z_2)+N(\phi^n|M\setminus M^{(n)}_{id})\leq
$$
$$
\leq \dim H_*(M^{(n)}_{id},\p{M^{(n)}_{id}};\Z_2)+N(\phi^n)\leq \dim H_*(M^{(n)}_{id},\p{M^{(n)}_{id}};\Z_2)+ \sum_{j} N((\phi)_j^n) +2k\leq Const
$$
by periodicity of $\phi_j$.
Taking the growth rate in $n$, we get
that asymptotic invariant  $ F^{\infty}(g)=1$.
\end{ex}
\begin{ex}
Let $\phi$ be an  hyperbolic automorphism  of 2-dimensional torus
defined by an integer matrix with eigenvalues $\lambda_1,\lambda_2, |\lambda_1|>1$. Then $ h(\phi)=log(|\lambda_1|)$. On the other hand
$N(\phi^n)=|det(I-A^n)|=|(1-\lambda_1^n)(1-\lambda_2^n)|$.
Hence  theorem \ref{thm:main1} implies
$$
 F^{\infty}(g) :=\grow(\dim\HF(\phi^n))= \limsup_{n \rightarrow \infty} |N(\phi^n)|^{1/n} =exp(h(\phi))=|\lambda_1| > 1
$$ 
\end{ex}

\subsection{Conjectures}

\subsubsection{Pseudo-Anosov mapping class}
 
\begin{conj}\label{conj:1}
For pseudo-Anosov  mapping class   $g \in \Gamma = \pi_0(Diff^+(M))$ we have
$$
\HF(g) =
\Z_2^{N (g)},  \,\,\,  \dim\HF(g)= N (g), \,\,\,  F^{\infty}(g)= \limsup_{n \rightarrow \infty} |N(g^n)|^{1/n}=h(\psi)=\lambda > 1
$$
where  $ N (g)$ denotes the Nielsen  number of  $g$ and $\psi$ is a standard(Thurston canonical form) representative of $g$,
i.e  there is  a monotone  representative $\phi\in\symp^m(M,\omega)$ of $g$
such that
$$
\HF(\phi) =
\Z_2^{N (\phi)},  \,\,\,  \dim\HF(\phi)= N (\phi), \,\,\, F^{\infty}(g) = \limsup_{n \rightarrow \infty} |N(\phi^n)|^{1/n}=h(\psi) =\lambda > 1
$$ 
\end{conj}

\begin{rk}

  For pseudo-Anosov  ``diffeomorphism'' we also have, as in theorems \ref{thm:main},\ref{thm:main0}
and
 \ref{thm:main1}, a topological separation of fixed points 
 \cite{Th, JG, I}, i.e the Nielsen number of pseudo-Anosov diffeomorphism
 equals to the number of
fixed points and there are  no connecting orbits between them.
But  pseudo-Anosov diffeomorphism is a symplectic
automorphism only on the complement of his fixed points set.
\end{rk}

Fathi and Shub \cite{flp} has proved the existence of Markov partitions for a
pseudo-Anosov homeomorphism $f$ representing mapping class $g$.The existence of Markov partitions implies that there is a symbolic dynamics for $(M ,f)$.This means that there is a finite set
$N$, a matrix $A=(a_{ij})_{(i,j)\in N\times N}$ with entries $0$ or $1$ and a surjective map $p:\Omega\rightarrow M$,where
$$
\Omega=\{(x_n)_{n\in \Z}: a_{x_nx_{n+1}}=1  ,  n\in \Z \}
$$
such that $p\circ \sigma =f\circ p$ where $\sigma$ is the shift (to the left) of the sequence $(x_n)$ of symbols.We have firstly:
$$
\# \fix\sigma ^n=\tr A^n.
$$
In general $p$ is not bijective.The non-injectivity of $p$ is due to the fact that the rectangles of the Markov partition can meet on their boundaries.To cancel the overcounting of periodic points on these boundaries,we use Manning's combinatorial arguments \cite{m} proposed in the case of Axiom A diffeomorphism . Namely, we construct finitely many subshifts of finite type ${\sigma _i}, i=0,1,..,m$, such that $\sigma_0=\sigma$, the other shifts semi-conjugate with restrictions of $f$ ,and signs $ \epsilon _i\in \{-1,1\}$ such that for each $n$
$$
\# \fix f^n=\sum_{i=0}^m \epsilon_i\cdot\#\fix \sigma_i^n =\sum_{i=0}^m \epsilon_i\cdot\tr A_i^n,
$$
where $A_i$ is transition matrix, corresponding to subshift of finite type $\sigma_i$.
For pseudo-Anosov homeomorphism of compact surface  $N(f^n)=\# \fix (f^n)$ for each $n>o$ \cite{Th}.So we have following trace formula for Nielsen numbers
$$
N(f^n)=\sum_{i=0}^m \epsilon_i\cdot\tr A_i^n.
$$

Conjecture \ref{conj:1} and this trace formula  immediately implie
the following
\begin{conj} For any  pseudo-Anosov mapping class $g$ the
symplectic zeta  function $F_g(z)$ is a rational function given
by the formula
\begin{equation}
F_g(z)=N_g(z)=N_f(z)=\prod_{i=0}^m\det(1-A_iz)^{-\epsilon_i}
\end{equation}

\end{conj}

\subsubsection{The general case.Concluding remarks}

\begin{conj}
For any  mapping class   $g \in \Gamma = \pi_0(Diff^+(M))$ there is
a monotone representative  $\phi\in\symp^m(M,\omega)$  with
respect to some $\phi$-invariant area form $\omega$  such that 
$$
\HF(\phi) =
H_*(M_{id},\p_{M_{id}};\Z_2)\oplus
\Z_2^{N(\phi|M\setminus M_{id})},
$$
where by $M_{id}$ we denote the union of the components
of $M\setminus\text{int}(U)$, where $\phi$ restricts to the identity.
Suppose $\psi$ is a standard( Thurston canonical form) representative
of $g$ and $\lambda $ is the largest stretching factor of
pseudo-Anosov pieses of $\psi$( $\lambda:=1 $ if there is no pseudo-Anosov
piece).Then asymptotic invariant
$$
 F^{\infty}(g):=\grow(\dim\HF(\phi^n)) =\lambda=h(\psi) = \limsup_{n \rightarrow \infty} |N(\psi^n)|^{1/n} 
$$ 
\end{conj}

\begin{rk}
(i)If $\phi \in\symp^m(M,\omega)$ has only non-degenerate fixed points,
then 
$$\# Fix(\phi) \geq \dim\HF(\phi)$$.
(ii) Due to P. Seidel  \cite{S1}   $\dim\HF(\phi)$  is a new symplectic invariant
of a four-dimensional symplectic  manifold  with nonzero first Betti number.
This 4-manifold produced from symplectomorphism $ \phi$ by a surgery construction which is a variation of earlier constructions due to McMullen-Taubes,
Fintushel-Stern and J. Smith.We hope that our  symplectic zeta functions 
and asymptotic invariant  also
give rise to a new  invariants of contact 3- manifolds
and symplectic 4-manifolds. 
\end{rk}


\begin{thebibliography}{10}
\bibitem{bbpt}
R.~Brooks, R.F.~Brown, J.~Pak, and D.H.~Taylor.
Nielsen numbers of maps of tori,
Proc. Am.Math.Soc., 52(1975), 398-400.
\bibitem{ds}
S.~Dostoglou and D.~Salamon.  Self dual instantons and holomorphic
  curves, Annals of Math., 139 (1994), 581--640.


\bibitem{EE}
C.~J. Earle and J.~Eells.
 The diffeomorphism group of a compact Riemann surface.
  Bull. Amer. Math. Soc., 73:557--559, 1967
\bibitem{flp}
A.~Fathi, F.~Laudenbach, and V.~Po{\'e}naru, Travaux de Thurston sur
  les surfaces, Ast{\'e}risque, vol. 66--67, Soc. Math. France, 1979.

\bibitem{f1}
A.L.~ Fel'shtyn.  Dynamical zeta functions, Nielsen theory and
Reidemeister torsion.
Memoirs of the American Mathematical Society, v.147, no.699, September 2000,
146 pages.
\bibitem{f2}
A.L.~Fel'shtyn.  Dynamical zeta functions and asymptotic expansions
in Nielsen theory.  Dynamical, Spectral, and Arithmetic Zeta
Functions. Contemporary Mathematics, v. 290, 2002, 67-81.
\bibitem{Floer1}
A.~Floer.
Morse theory for Lagrangian intersections.
J. Differential Geom., 28(1988), 513-547.

\bibitem{Floer}
A.~Floer.
 Symplectic fixed points and holomorphic spheres.
 Comm. Math. Phys., 120(2),575--611, 1989.

\bibitem{fri1}
D.~ Fried.
Lefschetz formula for flows,
The Lefschetz centennial conference.
Contemp. Math., 58 (1987), 19-69.

\bibitem{G}
R.~Gautschi.
Floer homology of algebraically finite mapping classes
 Preprint,April 2002. math.SG/0204032.

\bibitem{Gromov}
M.~Gromov.
Pseudoholomorphic curves in symplectic manifolds.
 Invent. Math.82(1985), 307-347.

\bibitem{J}
B.~Jiang.
 Fixed point classes from a differentiable viewpoint.
 In  Fixed point theory,  Lecture Notes in
  Math.,vol. 886, 163--170. Springer, 1981.
\bibitem{j1}
B. Jiang,
Estimation of the number of periodic orbits.
Pacific Jour. Math., 172(1996), 151-185.
\bibitem{JG}
B.~Jiang and J.~Guo.
 Fixed points of surface diffeomorphisms.
 Pac. J. Math., 160(1):67--89, 1993.
\bibitem{I}
N. V.~ Ivanov,
Nielsen numbers of maps of surfaces.
Journal  Sov. Math., 26, (1984).
 \bibitem{i1}
N. V. Ivanov,
Entropy and the Nielsen Numbers.
Dokl. Akad. Nauk SSSR 265 (2) (1982), 284-287 (in Russian);
English transl.:
Soviet Math. Dokl. 26 (1982), 63-66.
\bibitem {m}
A. Manning,
Axiom A diffeomorphisms have rational zeta function.
Bull. London Math. Soc. 3 (1971), 215-220.
 
\bibitem{MS}
D.~McDuff and D.~A. Salamon.
  Introduction to symplectic topology.
 Oxford Mathematical Monographs. Oxford Science Publications, 1998.
\bibitem{Mo}
J.~Moser.
 On the volume elements on a manifold.
  Trans. Amer. Math. Soc., 120:286--294, 1965.
\bibitem{Po}
M.~Po{\'z}niak.
 Floer homology, {N}ovikov rings and clean intersections.
 In  Northern California Symplectic Geometry Seminar, volume~196
of  Amer. Math. Soc. Transl. Ser.2, pages 119--181. American Mathematical
Society, 1999.
\bibitem{S}
P.~Seidel.
 Symplectic Floer homology and the mapping class group.
  Preprint, March 2001.
 math.SG/0010301.
\bibitem{S1}
P.~Seidel.
Braids and symplectic four-manifolds with abelian fundamental group.
Preprint, February 2002.
math. SG/0202135.
\bibitem{Th}
W.~P. Thurston.
 On the geometry and dynamics of diffeomorphisms of surfaces.
 Bull. Amer. Math. Soc., 19(2):417--431, 1988.


\end{thebibliography}
\end{document}